\documentclass[12pt]{amsart}
\usepackage{epsf}
\usepackage{amsmath}
\usepackage{amssymb}
\newtheorem{theo}{\bf Theorem}%[section]

\newtheorem{lemma}[theo]{\bf Lemma}
\newtheorem{conj}[theo]{\bf Conjecture}

\newtheorem{problem}[theo]{\bf Problem}
\newcommand{\DEF}{\buildrel {\rm def} \over =}

\begin{document}

\title[A two-sided secretary problem]
{Optimal stopping in a two-sided secretary problem}
\author{Kimmo Eriksson}
\address{IMA, M{\"a}lardalens h{\"o}gskola \\
   Box 883 \\
   SE-721 23 V{\"a}ster{\aa}s, Sweden}
\email{kimmo.eriksson@mdh.se}
\author{Jonas Sj{\"o}strand}
\address{Dept.~of Mathematics \\
   KTH \\
   SE-100 44 Stockholm, Sweden}
\email{jonass@kth.se}
\author{Pontus Strimling}
\address{IMA, M{\"a}lardalens h{\"o}gskola \\
   Box 883 \\
   SE-721 23 V{\"a}ster{\aa}s, Sweden}
\email{pontus.strimling@mdh.se} 
\date{July 1, 2004}

\begin{abstract}
In the ''secretary problem,'' well-known in the theory of optimal
stopping, an employer is about to interview a maximum of $N$
secretaries about which she has no prior information. Chow et
al.~proved that with an optimal strategy the expected rank of the
chosen secretary tends to $\prod_{k=1}^\infty
(1+\frac{2}{k})^{1/(k+1)} \approx 3.87$.

We study a two-sided game-theoretic version of this optimal
stopping problem, where men search for a woman to marry at the
same time as women search for a man to marry. We find that in the
unique subgame perfect equilibrium, the expected rank grows as
$\sqrt{N}$ and that, surprisingly, the leading coefficient is
exactly 1. We also discuss some possible variations.
\end{abstract}
\maketitle

\section{Introduction}
\noindent
A mathematical result of such general appeal that it sometimes
appears in ordinary newspapers (such as The Independent \cite{Co})
is the ''37 \% rule.'' This rule states that if you are walking
down a one-way street of restaurants of which you have no
knowledge beforehand and you have decided that, no matter what,
you will pick the hundredth restaurant if you haven't picked one
before that, then you should walk by 37 restaurants and then pick
the first one that you like better than all the previous
restaurants. This rule is optimal in the sense that you will
maximize your chance of picking the best of the hundred
restaurants.

\subsection{The secretary problem}
Although the restaurant setting may be the one most pertinent to
every-day life, the problem has become known as the
\emph{secretary problem} in the mathematics literature. Instead of
a hungry tourist judging restaurants, you are now an employer
interviewing candidates for a position as your secretary. After
each interview you must decide whether to hire this candidate or
continue the search process, in which case you cannot return to
this candidate. If for some reason you will look at $N$ candidates
at most, then you should never hire among the first $N/e$
candidates and then hire the first one that you like better than
all the previous ones. Since $1/e\approx 0.37$, this is the 37 \%
rule. The history of the secretary problem has been nicely told by
Ferguson \cite{Fe1}.

The secretary problem is the prime example of a question of
\emph{optimal stopping}. The theory of optimal stopping was
treated in a comprehensive way more than thirty years ago by Chow,
Robbins and Siegmund \cite{Ch1}, and more recently by Ferguson
\cite{Fe}. Our inspiration to the following piece of research came
from a short paper by Steven Finch on mathematical constants from
optimal stopping \cite{Fi}. After mentioning the secretary problem
and its constant $1/e$, Finch continues with the following much
less known variation, which is not treated in either of the
above-mentioned books: What is the optimal stopping strategy if
you want to minimize the expected rank $R_N$ of the chosen
secretary? (The $N$ candidates are ranked so that the best one has
rank 1, the second-best has rank 2, etc.) Lindley \cite{Li} found
an optimal strategy, and Chow et al.~\cite{Ch} proved that using
this strategy the expected rank $R_N$ tends to a constant as $N$
tends to infinity:
\[
\lim_{N\rightarrow \infty} R_N = \prod_{k=1}^\infty
\left(1+\frac{2}{k}\right)^{1/(k+1)} \approx 3.87
\]
We found this result doubly amazing --- not only does it defy
intuition that even with millions of candidates we can expect to
end up with the fourth-best just by using a stopping rule, but it
also seems unlikely that the constant could be expressed on such a
closed form.

\subsection{A two-sided secretary problem}
The article in The Independent on the 37 \% rule was based on an
interview with Peter Todd, a psychologist working on social
heuristics, i.e. behavioral rules of thumb. By laboratory
experiments as well as computer simulations, Todd and Miller
\cite{To, Du} have investigated simple heuristics for mate search
where both sexes have the possibility of saying no to a partner.
Mate search can therefore be described as a two-sided secretary
problem. There are several natural variants of the exact
specification. However, we have been unable to find any previous
mathematical treatment of any variant of the two-sided secretary
problem in the optimal stopping literature.

In this paper we are most interested in the following game
theoretic setting of the problem. There is a large universe of $U$
men and $U$ women. Every woman has a personal total preference
order on the men, and vice versa. In other words, for each woman
$w$ there is a permutation of the set of men representing her
personal rank order, and similarly for each man. The preferences
are unknown, only to be partially revealed through extensive
dating. Each of these persons are willing to search for a mate for
a maximum of $N$ rounds, where $N<<U$ (think of $N$ in the tens or
hundreds, and $U$ in the billions). In each round every woman will
date a new man, and they can either agree on marrying (and leave
the game) or at least one of them decides not to marry, in which
case they proceed to the next round and will never date each other
again. In the $N$th round, the players have no other option than
marrying. The condition that $N<<U$ will guarantee that all
players can count on new partners being available for all $N$
rounds, even though most will leave the game before the last
round.

Now, where do the permutations come from? We make the usual
assumption in the secretary problem that each permutation is drawn
randomly among all possible permutations. We also assume that a
random choice determines who will date who in each round (among
all unseen partners still available). This means that from the
viewpoint of any player entering round $r$, the rank of the next
date relative the $r-1$ partners already observed is a random
variable drawn from a uniform distribution on the set of ranks
from 1 to $r$. Let us call this the \emph{next-player principle}.

Finally, what is the goal of the players? It is not very
meaningful to maximize the probability of marrying the top-ranked
mate, since the player will only ever see a very small portion of
the universe, $N$ out of $U$. Instead, as in Lindley's variant of
the secretary problem \cite{Li,Ch} we assume that each player
wants to minimize the expected rank of the mate among the $N$
partners the player would meet if she completed all $N$ rounds.
Although the actual set of partners that a player would have met
is not determined if she does not play the whole game, the
expected rank can be computed by a simple formula that follows by
induction from the next-player principle:
\[
R_N(r,R_r) = \frac{N+1}{r+1}R_r,
\]
where $R_N$ is the expected rank in round $N$ of our current date
in round $r$, ranked $R_r$ among the $r$ partners observed up to
then. When several rank concepts are thrown around simultaneously,
we will refer to this as the \emph{$N$-rank}.

A \emph{strategy} in this two-sided secretary game is a rule that
says whether to propose marriage in round $r$ to a date of
observed rank $R_r$. Remember that the player cannot be certain
that the date will agree on marrying. (If players on the other
side always agree, then we would be back in the one-sided
secretary problem.)

\subsection{A cooperative problem}
To begin with, we can take the viewpoint of cooperative game
theory where players can make a binding agreement to play a
certain strategy.

\begin{problem} If the players can make a binding agreement
beforehand on a common strategy, what is the optimal choice? What is the
expected $N$-rank under this strategy?
\end{problem}

In this case the problem is a pure optimization problem of finding
the strategy on which to agree. Given every player's strategy, the
$N$-rank $R_N(r)$ a player can expect when she enters round $r$
can be computed by the following recurrence:
\begin{equation}\label{eq:fundamental}
R_N(r)=P[\mbox{marry}]\cdot \frac{N+1}{r+1}\cdot
E[R_r|\mbox{marry}] + (1-P[\mbox{marry}])\cdot R_N(r+1).
\end{equation}
The boundary condition is $R_N(N)=(N+1)/2$ since at the last round
each player can expect to obtain an average partner. As before,
$R_r$ denotes the rank of your current date among the $r$ partners
you have seen. We desire to compute $R_N(1)$, the expected
$N$-rank when we enter the game.

We have not solved Problem 1 completely, but we will motivate the
following conjecture which defines a new optimal stopping
constant:

\begin{conj}
In the cooperative two-sided secretary game, the asymptotics of
the expected $N$-rank when entering the game is given by
\[\lim_{N\rightarrow \infty}
\frac{R_N(1)}{\sqrt{N}} =\sqrt{27/32}\approx 0.92. \]
\end{conj}

\subsection{A noncooperative problem}
Now switch to the noncooperative  version of the game, where each
player tries to optimize her own outcome.

\begin{problem} What is the optimal strategy in the noncooperative two-sided
secretary game, and what is the expected $N$-rank under this
strategy?
\end{problem}

In noncooperative game theory, a strategy can only be optimal
given the strategies of the other players. If each player plays a
strategy that is optimal given the strategies of everybody else,
then we have a Nash equilibrium. A common problem in game theory
is that there exist very many Nash equilibria (such as ''always
marry in the first round''). The solution proposed by Selten \cite{Se}
is to assume that perfectly rational players will coordinate on a
subgame perfect equilibrium, that is, an equilibrium where all
players' strategies are optimal in every possible subgame. In the
two-sided secretary game, finding the unique subgame perfect
equilibrium is a recursive optimization problem: In the last
round, every player will accept marriage since it is always
preferable to not marrying at all; in the next to last round, a
player will want to marry if her current date is better than what
she would expect to obtain in the last round, etc.  We refer to
this as \emph{the} optimal strategy, which will be the same for
all players.

In this case we have been able to solve the problem completely.
The answer is surprisingly simple.

\begin{theo}\label{th:nash}
In the noncooperative two-sided secretary game, the asymptotics of
the expected $N$-rank when entering the game is given by
\[
\lim_{N\rightarrow \infty} R_N(1)/\sqrt{N} = 1.
\]
\end{theo}

This is our main result.

\subsection{The two-sided secretary problem with universal rank symmetry}
Due to the risk of rejections from good partners, each player in
the two-sided game is at a clear disadvantage compared to an
employer in a one-sided secretary problem.  This phenomenon is
inherent to two-sidedness, but it is reasonable to ask how small
the disadvantage can be made if preferences are dependent on each
other so that the probability of you liking me increases if I like
you.

An extreme case in this direction is \emph{universal rank
symmetry}. By this we mean that each woman $w$ and man $m$ have
the same universal rankings of each other, so that if $m$ is
ranked by $w$ as the $R$th best man in the universe, then $w$ is
ranked by $m$ as the $R$th best woman in the universe. (In
combinatorial terms, the women's permutations of the men
constitute a $U\times U$ Latin square since no two women can give
the same rank to the same man. The women's Latin square then
uniquely determines the men's Latin square.)

\begin{problem} What is the optimal strategy in the two-sided
secretary game with universal rank symmetry, and what is the
expected $N$-rank under this strategy?
\end{problem}

We conjecture that universal rank symmetry radically changes the
outcome, so that the expected $N$-rank tends to a small constant:

\begin{conj}
In the two-sided secretary game with universal rank symmetry, the
asymptotics of the expected $N$-rank when entering the game is
given by
\[\lim_{N\rightarrow \infty}
R_N(1) = {\rm constant} < 5. \]
\end{conj}

\subsection{Outline of paper}
We will work through the problems in the above order, starting
with our partial solution to the cooperative version. We will then
treat the noncooperative problem in the same way, but here we have
been able to compute the asymptotics using a method inspired by
Chow et al.~\cite{Ch}. After partially treating the case of
universal rank symmetry, we conclude with some discussion and open
problems.

\section{Problem 1: The optimal strategy under binding agreement}
\noindent
Let us assume that the players agree on a strategy, i.e.~a set of
thresholds $s_1$, $s_2$, \dots, $s_{N-1}$ together with the rule
that you must propose marriage to your date in round $r$ if and
only if his observed rank satisfies $R_r \le s_r$. Without loss of
generality we can always choose the threshold $s_r$ to be an
integer between 0 and $r$. Then, according to the next-player
principle, the probability that you will propose in round $r$ is
$s_r/r$. Hence, the probability of an agreement to marry is
\[
P[\mbox{marry}]=\left(\frac{s_r}{r}\right)^2,
\]
and the expected observed rank of your partner, given that you
propose, is
\[
E[R_r|\mbox{marry}]=\frac{s_r+1}{2}.
\]
We can plug these expressions into the fundamental recurrence
(\ref{eq:fundamental}). Setting $n=N-r$ (the number of rounds
remaining) and setting $\rho_n=R_N(N-n)\cdot 2/(N+1)$,  the
recurrence takes the shape
\begin{equation}\label{eq:fundamental-3}
\rho_n=\left(\frac{s_{N-n}}{N-n}\right)^2\cdot \frac{s_{N-n}+1}{N-n+1} +
\left(1-\left(\frac{s_{N-n}}{N-n}\right)^2\right)\cdot \rho_{n-1},
\end{equation}
with boundary condition $\rho_0=1$. We want to minimize
$\rho_{N-1}$. All factors and terms are nonnegative, so for each
$n$ from $1$ to $N-1$ we simply want to agree on the threshold
$s_{N-n}$ that will minimize $\rho_n$. To find this threshold,
differentiate the expression and find the zero to be
$s_{N-n}=\frac{2}{3}(-1+(N-n+1)\cdot\rho_{n-1})$. Then adjust the
threshold to the closest integer between 0 and $r$. For example,
setting $n=1$ we see that in the next-to-last round we shall
propose if our current date is among the best two thirds of all
partners we have seen.

For $n<<N$ we obtain the approximative recurrence
\[
\rho_n\approx \rho_{n-1}-\frac{4}{27}\rho_{n-1}^3, \qquad
\rho_0=1.
\]
Making the ansatz $\rho_n=A\cdot (n+B)^C$ and approximating the
difference $\rho_n-\rho_{n-1}$ by the derivative, we obtain the
solution $\rho_n=\sqrt{27/8}\cdot(n+4)^{-1/2}$.

\begin{figure}[htb]
\centerline{\epsfysize=100mm\epsfbox{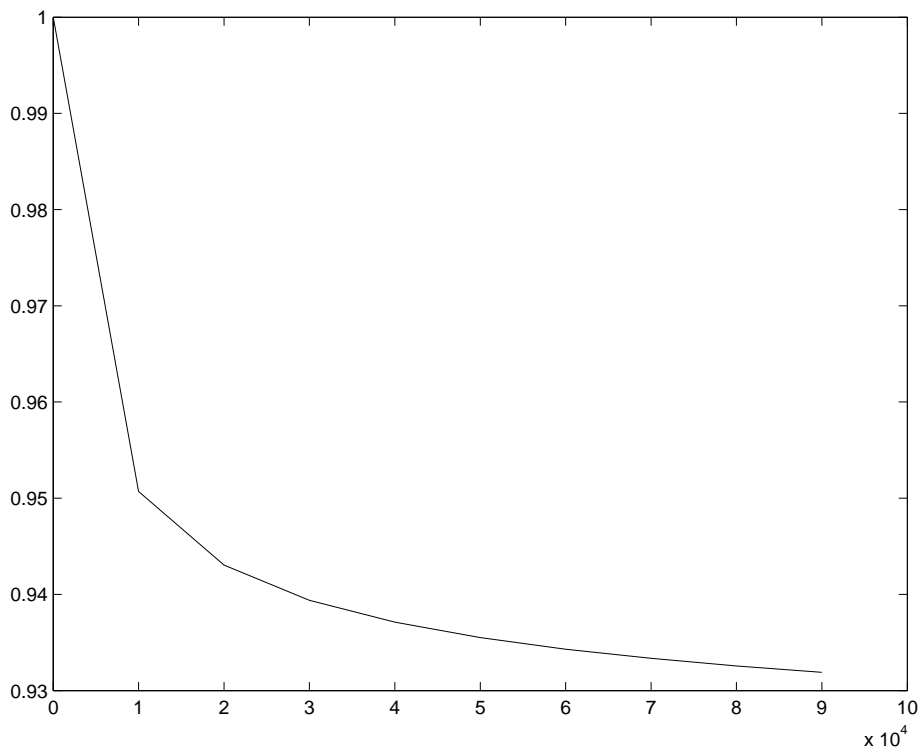}}
 \caption{\label{fig:agree}The graph of $R_N(1)/\sqrt{N}$ for the
two-sided secretary game under optimal binding agreement.}
\end{figure}

In the unlikely event that this crude approximation works all the
way for $n=0,1,\dots,N-1$, then a player's expected $N$-rank when
entering the game would be $R_N(1)=\frac{N+1}{2}\rho_{N-1}\approx
\sqrt{27N/32}$. Amazingly, computer calculations seem to vindicate
that $\lim_{N\rightarrow \infty} R_N(1)/\sqrt{N} \approx 0.92
\approx \sqrt{27/32}$, see the graph in Figure \ref{fig:agree}.

To summarize this section: We have found a recurrence giving the
optimal strategy under cooperation, and we conjecture that the
asymptotic behavior is given by $\lim_{N\rightarrow \infty}
R_N(1)/\sqrt{N} = \sqrt{27/32}$. It is possible that this problem
might yield to the same method that we use in the next section.

\section{Problem 2: The optimal noncooperative strategy}
\noindent
The strategy on which the players agreed in the previous section
is not a Nash equilibrium. This is obvious from the observation
that in the next-to-last round all players agreed to propose to
any partner among the best two thirds, while a rational player who
is not under any binding agreement would not accept a partner
worse than average, since she would expect an average partner in
the last round.

The thresholds that define the subgame perfect equilibrium
strategy in this noncooperative case are determined by what you
can expect to get if you decline to marry. You should accept
marriage in round $r$ if and only if the expected $N$-rank if you
marry now is less than or equal to the expected $N$-rank if you do
not marry:
\[
\frac{N+1}{r+1}\cdot R_r \le R_N(r+1),
\]
and the threshold $s_r$ should be the largest integer value of
$R_r$ satisfying the above inequality, so that
\[
s_r = \lfloor \frac{r+1}{N+1}\cdot R_N(r+1) \rfloor.
\]
We can assume all other players to reason in the same way, so that
we have $P[\mbox{marry}]=(s_r/r)^2$ as in the previous section.
Consequently, we can just plug this new value of $s_r$ into the
recurrence (\ref{eq:fundamental-3}). The same approximations now
yield
\[
\rho_n\approx \rho_{n-1}-\frac{1}{8}\rho_{n-1}^3, \qquad \rho_0=1,
\]
with an approximate solution of $\rho_n= 2(n+4)^{-1/2}$. As in the
cooperative case, the validity of this approximation is supported
by computer calculations indicating that
\[
\lim_{N\rightarrow \infty} R_N(1)/\sqrt{N} = 1.
\]
This is our Theorem~\ref{th:nash}. We will prove this result
using the same approach as Chow et al.~\cite{Ch}, although we will
encounter different technical difficulties than they did.

\subsection{Proof of Theorem~\ref{th:nash}}
For convenience, introduce $i=r-1$ (so that the range is
$i=0,1,\dots,N-1$) and $c_i = R_N(i+1)$. We wish to prove that
$c_0/\sqrt{N}\rightarrow 1$. Define the unrounded threshold
\[
t_i = \frac{i+1}{N+1}\cdot c_i,
\]
so that $s_i=\lfloor t_i \rfloor$. Then recurrence
(\ref{eq:fundamental}) can be rewritten as
\begin{equation}\label{eq:st}
t_{i-1}=\frac{s_i^2(s_i+1)+2(i^2-s_i^2)t_i}{2i(i+1)}; \qquad
t_{N-1}=\frac{N}{2}.
\end{equation}
For any $i$, define $\alpha_i=t_i-s_i$. Then $0\le \alpha_i < 1$,
and (\ref{eq:st}) becomes
\begin{equation}\label{eq:t}
t_{i-1}=\frac{(t_i-\alpha_i)^2(t_i-\alpha_i+1)+2(i^2-(t_i-\alpha_i)^2)t_i}
{2i(i+1)}.
\end{equation}

It is trivial to see that
$a_it_i+(1-\alpha_i)[t_i^2+\alpha_i(t_i-\alpha_i)]\ge 0$. If we add
this number to the nominator of (\ref{eq:t}) we obtain the upper
bound
\begin{equation}\label{eq:ub}
t_{i-1}\le T_i(t_i) \DEF
\frac{-t_i^3+2t_i^2+2i^2t_i}{2i(i+1)}.
\end{equation}
Similarly, if we subtract
$\alpha_i[(t_i-1)^2+\alpha_i(t_i-\alpha_i)]+(t_i-\alpha_i)+\alpha_i^2\ge
0$ from the nominator of (\ref{eq:t}) we obtain the lower bound
\begin{equation}\label{eq:lb}
t_{i-1}\ge \tau_i(t_i) \DEF
\frac{-t_i^3+t_i^2+(2i^2-1)t_i}{2i(i+1)}.
\end{equation}
It is easy to show that, if $i\geq2$, the functions $T_i(t)$
and $\tau_i(t)$ are
both increasing for all $t$ in the range $[0,\sqrt{2/3}i]$, and
that this interval includes all values $t$ can possibly attain, so
we can use these recursive inequalities in order to prove explicit
lower and upper bounds.

\begin{lemma}\label{lm:ub}
Let $f(N)$ be any real function with
$\lim_{N\rightarrow\infty}f(N)=\infty$.
For sufficiently large $N$ the upper bound
\[
 t_i \le \frac{i+\sqrt{i}}{\sqrt{N-i+3}},
\]
holds for all $i$ in the interval $f(N)\leq i\leq N-1$.
\end{lemma}
\begin{proof}
Proof by backwards induction. The lemma holds for $i=N-1$, since
\[
t_{N-1}=N/2\le \frac{N-1+\sqrt{N-1}}{\sqrt{4}}
\]
for all $N\ge 2$. Assuming the lemma holds for a given $i$ we have
\[
t_{i-1} \le T_i(t_i) \le
T_i\big(\frac{i+\sqrt{i}}{\sqrt{N-i+3}}\big),
\]
thanks to (\ref{eq:ub}) and the fact that $T_i$ is an increasing
function.
Hence, to conclude the induction step we only need to
prove that
\[
T_i\big(\frac{i+\sqrt{i}}{\sqrt{N-i+3}}\big) \le
\frac{i-1+\sqrt{i-1}}{\sqrt{N-i+4}}.
\]
This can be verified through several steps of computations in
Maple (see Appendix).
\end{proof}

The lower bound is trickier. We have not been able to find a lower
bound for which the induction step works, but by adding a constant
term of 0.148 to the denominator we obtain a lower bound that
works from $N-22$ and downwards.

\begin{lemma}\label{lm:lb}
For sufficiently large $N$, the lower bound
\[
t_i\ge \frac{i+1}{\sqrt{N-i+3}+0.148}
\]
holds for all $i$ in the interval $\sqrt{N}+1\leq i\leq N-22$.
\end{lemma}
\begin{proof}
For large $N$ we can compute an approximation of $t_{N-22}$ by
using the recurrence 21 steps, keeping only the most significant term
in each step: If $t_{N-k}=a_kN+o(N)$ we have
$s_{N-k}=a_kN+o(N)$ too, so our recurrence~(\ref{eq:st})
gives $t_{N-k-1}=\frac{-a_k^3+2a_k}{2}N+o(N)$.
Maple gives
\[
t_{N-22}\approx 0.19427N \ge 0.19425N - 4.07925 =
\frac{(N-22)+1}{\sqrt{22+3}+0.148}.
\]
Again we proceed by backwards induction. Assuming the lemma holds
for a given $i$ we have
\[
t_{i-1} \ge \tau_i(t_i) \ge
\tau_i\left(\frac{i+1}{\sqrt{N-i+3}+0.148}\right),
\]
thanks to (\ref{eq:lb}) and the fact that $\tau_i$ is an
increasing function. It remains for us to prove that
\[
\tau_i\left(\frac{i+1}{\sqrt{N-i+3}+0.148}\right) \ge
\frac{i}{\sqrt{N-i+4}+0.148}
\]
for large $N$. Again, several steps in Maple verifies this
inequality, see Appendix.
\end{proof}

We are now in a position to finish the proof of
Theorem~\ref{th:nash}.
In the case $t_i<1$ we always have $s_i=0$ so no partners are ever
accepted and the expected $N$-rank $c_i$ remains constant down to
$c_0$. Let $i_{\rm crit}$ be the greatest $i$ with $t_i<1$, so that
$c_0=c_1=\cdots=c_{i_{\rm crit}}$.

Lemma~\ref{lm:lb} gives that, for large $N$,
$$t_{\lceil\sqrt{N}+1\rceil}\geq
\frac{\lceil\sqrt{N}+1\rceil+1}{\sqrt{N-\lceil\sqrt{N}+1\rceil+3}+0.148}>1$$
so $i_{\rm crit}<\sqrt{N}+1$.
Lemma~\ref{lm:ub} with $f(N)=N^{1/2}-N^{1/3}$ gives, for large $N$,
$$t_{\lceil N^{1/2}-N^{1/3}\rceil}\leq
\frac{\lceil N^{1/2}-N^{1/3}\rceil+\sqrt{\lceil N^{1/2}-N^{1/3}\rceil}}
{\sqrt{N-\lceil N^{1/2}-N^{1/3}\rceil+3}}$$
$$=\sqrt{\frac{N-2N^{5/6}+o(N^{5/6})}{N-N^{1/2}+o(N^{1/2})}}<1$$
so $i_{\rm crit}+1\geq N^{1/2}-N^{1/3}$.

Since $N^{1/2}-N^{1/3}\leq i_{\rm crit}+1<\sqrt{N}+2$,
Lemma~\ref{lm:ub} with $f(N)=N^{1/2}-N^{1/3}$ gives that
$$1\leq t_{i_{\rm crit}+1}\leq
\frac{i_{\rm crit}+1+\sqrt{i_{\rm crit}+1}}{\sqrt{N-i_{\rm crit}+2}}
\leq\frac{\sqrt{N}+2+\sqrt{\sqrt{N}+2}}{\sqrt{N-\sqrt{N}+1}}\rightarrow1$$
as $N\rightarrow\infty$.
Then $s_{i_{\rm crit}+1}=1$, and
by~(\ref{eq:st}) we get $t_{i_{\rm crit}}\rightarrow1$.

Recall that, by definition, $t_i = c_i(i+1)/(N+1)$ so that $c_{i_{\rm
crit}}=t_{i_{\rm crit}}(N+1)/(i_{\rm crit}+1)$. Thus, we have
\[
\frac{c_0}{\sqrt{N}} = \frac{c_{i_{\rm crit}}}{\sqrt{N}}
=\frac{N+1}{i_{\rm crit}+1}\cdot \frac{t_{i_{\rm crit}}}{\sqrt{N}}
\rightarrow 1
\]
since $t_{i_{\rm crit}}\rightarrow1$ and
$N^{1/2}-N^{1/3}-1\leq i_{\rm crit}<\sqrt{N}+1$.

\subsection{A social dilemma}
This game illustrates the game-theoretic concept of a \emph{social
dilemma}. All players would like to agree on the cooperative
strategy, where you accept to marry quite a lot of partners for
the good of the group. However, when an individual player finds
herself in a position where the strategy calls on her to marry a
date that she finds below her expectations, she is tempted to
reject this partner and optimize her own good instead. But if all
players do that, then the expected outcome is worse for all of
them; for large $N$ the expected change for the worse is about 8
percent.

\section{Problem 3: The optimal strategy under universal rank symmetry}
\noindent
If the preferences of all players satisfy universal rank symmetry,
then it should be easier for the players to find mutually
acceptable agreements.

As in the previous problem, a rational player will accept a
partner in round $r$ if the observed rank does not exceed the
threshold $s_r = \lfloor \frac{r+1}{N+1}\cdot R_N(r+1) \rfloor$.
The new circumstance in the current setting is that the events
''your observed rank of me is at most $s_r$'' and ''my observed
rank of you is at most $s_r$'' are no longer independent. Hence
the probability $P[\mbox{marry}]=P[\mbox{both ranks } \le s_r]$
must be found by other means. Integration over the unknown global
rank, which is the same for both players, yields
\begin{equation}\label{eq:sur-P}
P[\mbox{marry}]=\sum_{k=0}^{s_r-1}\sum_{\ell=0}^{s_r-1}
\frac{\binom{r-1}{k}\binom{r-1}{\ell}}{\binom{2(r-1)}{k+\ell}(2r-1)}
\end{equation}
Similarly, the expected rank of a partner given the fact that you
both agree to marry is no longer given by the arithmetic mean in
the rank interval $[1,s_r]$ but is more favorable:
\begin{equation}\label{eq:sur-E}
E[R_r|\mbox{marry}]=\frac{r}{s_r}\sum_{k=0}^{s_r-1}(k+1)
\binom{r-1}{k}\sum_{\ell=0}^{s_r-1}
\frac{\binom{r-1}{\ell}}{\binom{2(r-1)}{k+\ell}(2r-1)}
\end{equation}
The recurrence resulting from plugging these expressions into
(\ref{eq:fundamental}) determines the optimal strategy under
universal rank symmetry.

We have not yet been able to find the limit of the expected
$N$-rank in this case, but computer calculations indicate that it
approaches a small constant (less than 5). Thus, making the
preferences dependent of each other in this way seems to have
changed the behavior of the expected rank so that it resembles the
behavior in the one-sided problem.

\section{Discussion and open problems}
\noindent
For now, we have to leave our conjectures of the asymptotic
behavior of the expected $N$-rank in Problem 1 and 3 as open
problems. Another open problem is to find the optimal strategy
when the universe is small, $U=N$, so that you do not know during
a date whether there will be any more dates or if all possible
partners will already be married if you reject the one you are
currently entertaining. In this case, another reasonable version
is obtained if we lift the restriction on dating the same person
twice so that in each round one would simply select a random date
among all remaining unmarried players of the opposite sex.

In the context of the two-sided secretary problem, it is
reasonable to briefly discuss the theory of \emph{stable
matching}, pioneered by Gale and Shapley \cite{Ga} and later
comprehensively treated by Roth and Sotomayor \cite{Ro}. Also this
theory studies men and women who have preferences on each other.
The difference is that the stable matching theory assumes that
each player knows all her preferences in advance, while in the
secretary problem the preferences are only revealed slowly as the
player meets new potential partners. On the other hand, in the
secretary problem we have no guarantees that the marriages will be
stable when new partners are met later in life. It would be
interesting to study some measure of how stable these marriages
will be under various conditions.

In a recent study, Caldarelli and Capocci \cite{Ca} have studied
stable marriages under the assumption that people's preferences
are influenced by a commonly appreciated trait such as beauty.
They find that this condition favors the very beautiful players
while all others are worse off than in a world where preferences
are random. Universal rank symmetry is an assumption in the
opposite direction. True preferences are likely to reflect a
mixture of randomness, a common sense of beauty, and an
''I-like-you-if-you-like-me'' component like our universal rank
symmetry. Modelling such sets of permutations is a challenge to
combinatorialists.

\section*{\bf Appendix}
\noindent
Here we prove two lemmas that are needed in the proof
of Theorem~\ref{th:nash}.
\begin{lemma}
Let $f(N)$ be any real function with
$\lim_{N\rightarrow\infty}f(N)=\infty$. Then,
for sufficiently large $N$, the inequality
$$
T_i\big(\frac{i+\sqrt{i}}{\sqrt{N-i+3}}\big) \le
\frac{i-1+\sqrt{i-1}}{\sqrt{N-i+4}}
$$
holds for all $f(N)\leq i\leq N-1$.
\end{lemma}

\begin{proof}
After the substitution $z=\sqrt{N-i+3}$ our inequality transforms to
\begin{equation}\label{eq:olik}
T_i\big(\frac{i+\sqrt{i}}{z}\big) \le
\frac{i-1+\sqrt{i-1}}{\sqrt{z^2+1}}
\end{equation}
and the condition $i\leq N-1$ transforms to $z\geq 2$. We will
prove that the inequality~(\ref{eq:olik}) holds for $z\geq1$ and
sufficiently large $i$. A large $N$ implies a large $i$ since
$f(N)\leq i$, so the lemma will follow.

Evaluating $T_i$ in~(\ref{eq:olik}) yields
$$g(i,z)\DEF \frac{i-1+\sqrt{i-1}}{\sqrt{z^2+1}}
- \frac{-\frac{(i+\sqrt{i})^3}{z^3}+\frac{2(i+\sqrt{i})^2}{z^2}
+\frac{2i^2(i+\sqrt{i})}{z}}{2i(i+1)} > 0.$$
We see that
$$g(i,1)=
(\sqrt{2}-1)i^3 + i^{5/2} + (\sqrt{2(i-1)}+1)(i^2+i)
-3i^{3/2} - (2+\sqrt{2})i$$
which is positive for large $i$.
Since $g(i,z)$ is continuous it suffices to show that
$g(i,z)\neq0$ for $z\geq1$ and $i$ large. The zeros of $g(i,z)$ are
the same as the zeros of
$$\left(\frac{i-1+\sqrt{i-1}}{\sqrt{z^2+1}}\right)^2
- \left(\frac{-\frac{(i+\sqrt{i})^3}{z^3}+\frac{2(i+\sqrt{i})^2}{z^2}
+\frac{2i^2(i+\sqrt{i})}{z}}{2i(i+1)}\right)^2.
$$
We multiply this by $(2i(i+1))^2(z^2+1)z^6$ and obtain a
polynomial of degree 6 in $z$:
\begin{eqnarray*}
q(z) & \DEF & 4(2(i^5+i^4-i^3-i^2)\sqrt{i-1}-2i^{11/2}-i^4-i^3)z^6 \\
& - & 8(i^5+3i^{9/2}+3i^4+i^{7/2})z^5 \\
& + & 4(2i^{11/2}+5i^5+4i^{9/2}-4i^{7/2}-6i^3-4i^{5/2}-i^2)z^4 \\
& + & 4(-i^5-i^{9/2}+4i^4+8i^{7/2}+5i^3+i^{5/2})z^3 \\
& + &(3i^6+10i^{11/2}+9i^5-4i^{9/2}-15i^4-22i^{7/2}-25i^3-16i^{5/2}-4i^2)z^2\\
& + & 4(i^5+5i^{9/2}+10i^4+10i^{7/2}+5i^3+i^{5/2})z \\
& - & (i^6+6i^{11/2}+15i^5+20i^{9/2}+15i^4+6i^{7/2}+i^3)
\end{eqnarray*}
We must show that $q(z)\neq0$ for $z\geq1$ and $i$ large. Using a computer
program like Maple this is just a matter of verification:
\begin{itemize}
\item
Compute the fourth derivative $q^{(4)}(z) = a_2(i)z^2+a_1(i)z+a_0(i)$.
\item
Verify that
$\left(\frac{a_1(i)}{2a_2(i)}\right)^2-\frac{a_0(i)}{a_2(i)} < 0$
for large $i$. This implies that $q^{(4)}(z)$ has no real roots.
\item
Check that $q^{(4)}(0)>0$ for large $i$ so that
$q^{(4)}(z)>0$ everywhere. This is immediately evident from
the expression for $q(z)$ above.
\item
Check that $q'''(1)$, $q''(1)$, $q'(1)$, and $q(1)$ are
positive for large $i$.
\end{itemize}
This shows that $q(z)>0$ for all $z\geq1$ for sufficiently large $i$.
\end{proof}

\begin{lemma}
For sufficiently large $N$, the inequality
$$
\tau_i\big(\frac{i+1}{\sqrt{N-i+3}+0.148}\big) \ge
\frac{i}{\sqrt{N-i+4}+0.148}
$$
holds for all $i$ in the interval $\sqrt{N}+1\leq i\leq N-22$.
\end{lemma}

\begin{proof}
Put $\varepsilon=0.148$. After the substitution $z=\sqrt{N-i+3}+\varepsilon$
our inequality transforms to
\begin{equation}\label{eq:grundolikhet}
\tau_i\big(\frac{i+1}{z}\big)\ge
\frac{i}{\sqrt{(z-\varepsilon)^2+1}+\varepsilon}.
\end{equation}
The interval $\sqrt{N}+1\leq i\leq N-22$ transforms to
\begin{equation}\label{eq:intervall}
5\leq z-\varepsilon\leq\sqrt{(i-1)^2-i+3}.
\end{equation}
Evaluating $\tau_i$ in~(\ref{eq:grundolikhet}) and multiplying by
$2iz^3(\sqrt{(z-\varepsilon)^2+1}+\varepsilon)>0$ yields
\begin{eqnarray*}
p(i) & \DEF &
((2z^2-1)\sqrt{(z-\varepsilon)^2+1}-(2z^3+(1-2z^2)\varepsilon))i^2 \\
& + & ((z-2)(\sqrt{(z-\varepsilon)^2+1}+\varepsilon)i \\
& - & (z^2-z+1)(\sqrt{(z-\varepsilon)^2+1}+\varepsilon)\\
& \geq & 0,
\end{eqnarray*}
a quadratic polynomial in $i$.

Let us first show that the leading
coefficient of $p(i)$ is positive, i.\ e.
$$
(2z^2-1)\sqrt{(z-\varepsilon)^2+1} > 2z^3+(1-2z^2)\varepsilon.
$$
After squaring, expanding, and collecting we get $4\varepsilon
z^3-3z^2-2\varepsilon z+1>0$. The zeros of this third-degree
polynomial in $z$ are $z_1\approx -0.592$, $z_2\approx 0.559$, and
$z_3\approx 5.100$, all of which are less than $5+\varepsilon$.

There are always real roots to $p(i)$ since
its constant term is negative for $z\geq5+\varepsilon$.
For each $z$, let $i(z)$ be the greater of the two roots of $p(i)$.
We can, of course, write down
an explicit formula for $i(z)$ and with a computer program like
Maple it is easy to check that
$\lim_{z\rightarrow\infty}i(z)-z=1-\varepsilon$. Thus, for any $\delta>0$
there is a $Z_\delta$ such that
$p(i)>0$ for $i>z+1-\varepsilon+\delta$ for all $z>Z_\delta$.
We choose $\delta=1/4$. In the interesting interval~(\ref{eq:intervall})
we have $z+1-\varepsilon+1/4\leq\sqrt{(i-1)^2-i+3}+3/4<i$, where the last
inequality follows from $(i-1)^2-i+3<(i-3/4)^2$ which is true if $i\geq3$.

In the interval $Z_{1/4}<z\leq\varepsilon+\sqrt{(i-1)^2-i+3}$ we know
that $p(i)>0$ so we only have to worry about the interval
$5+\varepsilon\leq z\leq Z_{1/4}$.
Since $i(z)$ is a continuous function we can define
$I=\sup\{i(z)\ :\ 5+\varepsilon\leq z\leq Z_{1/4}\}$. Then
$p(i)>0$ in the interval~(\ref{eq:intervall}) for all $i>I$. A large $N$
implies a large $i$ since $\sqrt{N}+1\leq i$, so the lemma follows.
\end{proof}

\section{Acknowledgment}
\noindent
Research partially supported by the European Commission's IHRP
Programme, grant HPRN-CT-2001-00272, ''Algebraic Combinatorics in
Europe.''

\end{document}